%% file: float_comment_arXiv.tex
\documentclass[12pt]{amsart}

\usepackage{amssymb}

\usepackage{epsfig,color}

\newtheorem{defin}{Definition}

\newtheorem{rem}{Remark}
\newtheorem{quest}{Question}

\newcommand{\RR}{{\mathbb R}^2}

\newcommand{\bo}{\partial} 

\newcommand{\const}{\mbox{const}} 

\newcommand{\al}{\alpha}

\newcommand{\ga}{\gamma}\newcommand{\Ga}{\Gamma}
\newcommand{\de}{\delta}

\newcommand{\Om}{\Omega}

\newcommand{\het}{\theta}

\begin{document}

\bibliographystyle{plain}

\title[Addendum to capillary floating]
{Addendum to: Capillary floating and the billiard ball problem}

\author{Eugene Gutkin}

\address{Department of Mathematics, Nicolaus Copernicus University, Chopina 12/18, Torun 87-100, Poland;
Institute of Mathematics of Polish Academy of Sciences,
Sniadeckich 8, Warsaw 00-956, Poland}
\email{gutkin@mat.umk.pl,gutkin@impan.pl}



\date{\today}

\begin{abstract}
We compare the results of our earlier paper on the floating in
neutral equilibrium at arbitrary orientation in the sense of
Finn-Young with the literature on its counterpart in the sense of
Archimedes. We add a few remarks of personal and social-historical
character.
\end{abstract}

\maketitle


%
Humans were experimenting with floating long before they were able
to come up with any understanding of this phenomenon. Current
theory as it has developed over the past centuries is quite
complicated, and various simplified models have been proposed.
Each such model singles out some of the physical aspects of
floating and disregards others. For instance, the ``archimedean''
floating accounts for buoyancy in a gravity field but disregards
the fluid surface tension. On the other hand, Thomas
Young~~\cite{Young}\footnote{Of the Rosetta stone fame.} based the
theory on surface tension and concluded that the {\em contact
angle} of the liquid surface with the body must be
prescribed.\footnote{See also a pioneering study of capillary
floating by Laplace \cite{La1806}.} We refer the reader to the
papers of R. Finn et al \cite{F09i,F09ii,F11} for elaborations.
Here we are concerned only with the geometric aspects of floating
models.


\medskip

A case of particular interest appears when the floating solid is
an infinite cylinder resting horizontally on the liquid. A
three-dimensional floating model then reduces, by translational
invariance, to a two-dimensional theory, expressed by suitable
conditions on the cross-section of the cylinder, $\Om\subset\RR$,
a bounded domain with a piecewise smooth boundary $\bo\Om$. We
identify the set of directions of oriented straight lines in $\RR$
with the circle $\{\al:0\le\al<2\pi\}$. Finn introduced the {\em
neutral equilibrium} model in which the fluid is horizontal and
meets the body at a prescribed angle, say $\pi-\ga$. See
figure~~\ref{finn_young_float_fig}. This floating model takes into
account the capillary forces and the liquid surface tension but
disregards the gravity. We will simply refer to it as the {\em
Finn-Young floating}. According to this model, $\Om$ {\em floats
in every orientation at the contact angle $\pi-\ga$} if and only
if every directed chord making angle $\ga$ with $\bo\Om$ at one
end, makes the same angle with $\bo\Om$ at the other end. In
\cite{F09i} Finn posed
\begin{quest}      \label{finn_quest}
What convex domains $\Om$ (other than the round disc) satisfy this
condition and for what $0<\ga<\pi/2$?
\end{quest}
%


Using the unpublished work \cite{Gut93}, the author provided a
fair amount of information on this in \cite{Gut12}: There is a
dense countable set $\Ga\subset(0,\pi/2)$, and for any $\ga\in\Ga$
there is a real analytic one-parameter family of distinct, smooth,
strictly convex domains $\Om_{\ga,\tau}$ that float in every
orientation at the contact angle $\pi-\ga$. The set $\Ga$ is as
follows:
\begin{equation}     \label{trig_eq}
\Ga=\cup_{n=2}^{\infty}\{0<\ga<\pi/2:\tan n\ga\,=\,n\tan\ga\}.
\end{equation}
The domains $\Om_{\ga,\tau}$ are explicitly described via the
Fourier coefficients of their radius of curvature functions
\cite{Gut12}. Incidentally, these results also bear on the {\em
billiard ball problem}, which yielded basic open questions in
geometry and analysis \cite{Gut03}.

\begin{figure}[htbp]
\input{finn_young_float.pstex_t}
\begin{center}
\caption{Finn-Young floating at the orientation $\al$ with the
contact angle $\pi-\ga$.} \label{finn_young_float_fig}
\end{center}
\end{figure}
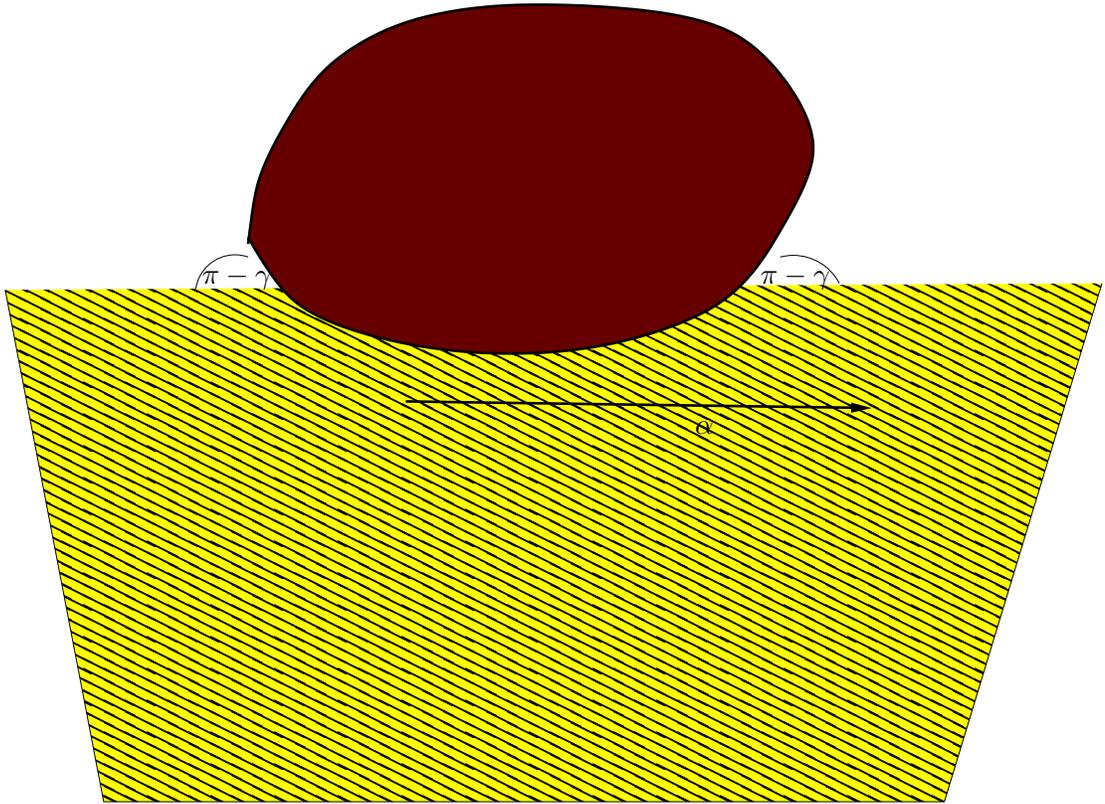
\begin{rem}     \label{four_rem}
{\em The work~~\cite{RMBC} contains an earlier surprising result
on the Finn-Young floating. It says that every smooth, strictly
convex plane domain $\Om$ floats in neutral equilibrium at any
contact angle in at least four distinct orientations. The bound is
the best possible. }
\end{rem}

\medskip

We will now turn to the archimedean floating. If the fluid surface
is assumed flat, the round ball floats in every orientation and at
every density less than the density of the liquid. Suppose that a
bounded solid floats in every orientation at a fixed density. One
of the problems of the famous ``Scottish Caf\'e'' \cite{Scot},
attributed to S. Ulam, asks if the solid is a round ball. The $2D$
version of this problem, also attributed to Ulam~~\cite{Au38},
asks the same question for bounded plane domains.\footnote{The
two-dimensional version comes from the archimedean floating of
infinite cylinders.} For reader's convenience, we will now
formulate the corresponding geometric problem.

Let $\Om\subset\RR$ be a bounded region with a piecewise smooth
boundary $\bo\Om$. The following sets are in natural one-to-one
correspondences: Ordered pairs of points in $\bo\Om$, arc segments
in $\bo\Om$, and directed chords in $\Om$ with the endpoints on
$\bo\Om$. We will use the notation $|\ |$ for the arc length; in
particular, $|\bo\Om|$ is the perimeter of $\Om$. Let
$\Om(A)\subset\RR$ be the region enclosed between the arc segment
$A$ and the corresponding chord $C(A)$.

\begin{defin}   \label{archim_float_def}
Let $0<\de\le 1/2$. We say that $\Om$ {\em floats (in the
archimedean sense) in equilibrium in every orientation at the
density $\de$} if for all arc segments $A\subset\bo\Om$ satisfying
$|A|=\de|\bo\Om|$ the areas $|\Om(A)|$ are equal.
\end{defin}
The archimedean counterpart of Question~~\ref{finn_quest} is thus
the following.
\begin{quest}      \label{archimed_que}
What plane domains (not necessarily  convex), other than the round
disc, satisfy the condition of Definition~~\ref{archim_float_def},
and what $\de$ are possible?
\end{quest}
%

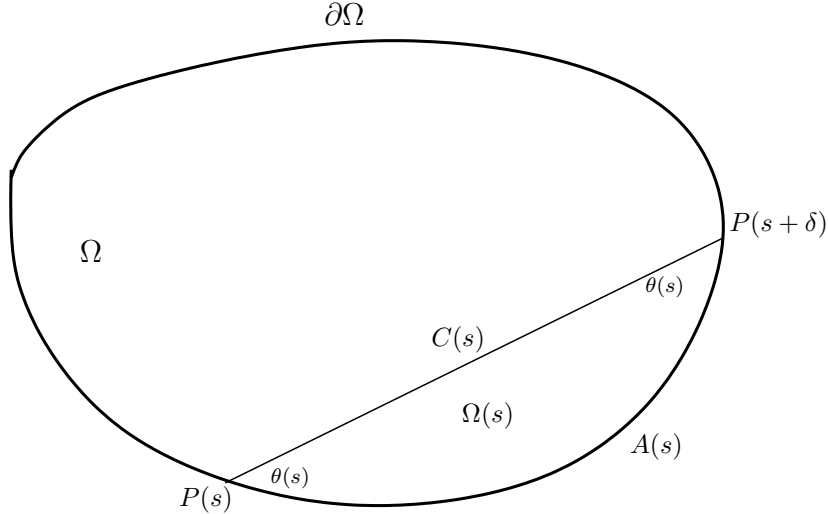
\begin{figure}[htbp]
\input{arch_float.pstex_t}
\begin{center}
\caption{Archimedean floating in every orientation at density
$\de$.} \label{arch_float_fig}
\end{center}
\end{figure}

Question~~\ref{archimed_que} has been studied much more than
Question~~\ref{finn_quest}, but despite the efforts of geometers
that go back to the 1920s
\cite{Sa34,Ge36,Au38,ZaKo_french,Ru_french,Ta06,Va09,Zi20}, the
results are less complete. We will now very briefly discuss the
literature on Question~~\ref{archimed_que} using the term {\em
floating domains} for plane regions satisfying the conditions of
Definition~~\ref{archim_float_def}. See
figure~~\ref{arch_float_fig} for the  notation. We assume without
loss of generality that $|\bo\Om|=1$, and parameterize $\bo\Om$ by
the arc length $0\le s < 1$. Let $P(s),P(s+\de)$ be the endpoints
of the moving arc segment $A(s)$. By elementary differential
geometry, the condition $|\Om(s)|=\const$ holds if and only if
$|C(s)|=\const$ if and only if the angles between $C(s)$ and
$\bo\Om$ at the points $P(s),P(s+\de)$ are equal.\footnote{See,
e.g., equation (1.1) in~~\cite{GK95}.} Let $\het(s)$ be the angle
at $P(s)$. Note that $\het(\cdot)\ne\const$, in general. Moreover,
it follows from \cite{Kov80} that $\het(\cdot)=\const\ne \pi/2$
implies that $\Om$ is a disc.

Question~~\ref{archimed_que} has been fully answered for convex
$\Om$ and $\de=1/2$. In~~\cite{Au38} H. Auerbach describes the
curves $\bo\Om$ satisfying this condition; he calls them {\em the
Zindler curves}~~\cite{Zi20}. Auerbach characterizes these curves
via the Fourier expansion of their radius of curvature functions.
He emphasizes an analogy between the Zindler curves and the {\em
curves of constant width}.

\begin{figure}[htbp]
\input{za_ko.pstex_t}
\begin{center}
\caption{The non-convex floating domain corresponding to the
equilateral triangle.} \label{ZaKo_fig}
\end{center}
\end{figure}
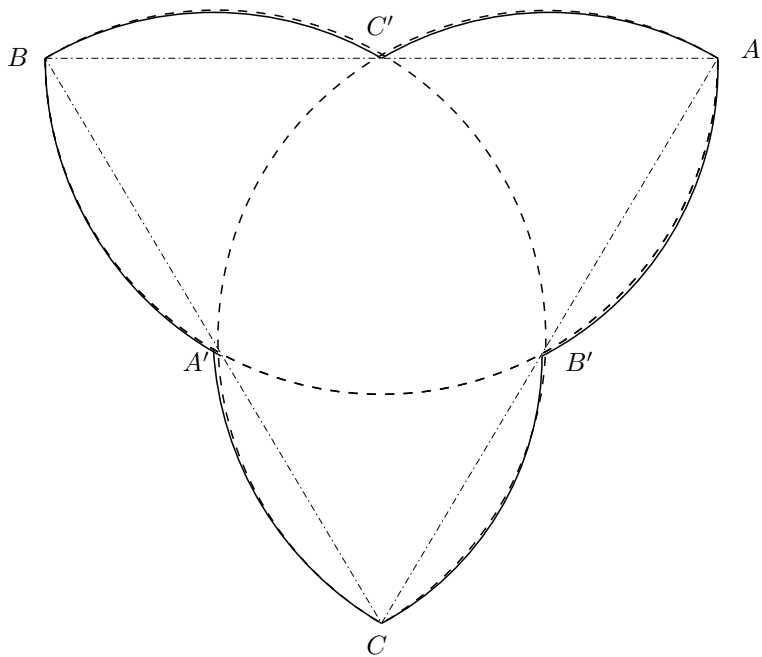

The work~~\cite{ZaKo_french} contains a surprising construction of
a class of nonconvex domains satisfying the conditions of
Question~~\ref{archimed_que} with $\de=1/2$. Let $P=A_1\cdots
A_{2n}\subset\RR$ be a $2n$-gon such that i) the sides of $P$ have
the same length; ii) there is $1<k\le n$ such that all diagonals
$A_iA_{i+k}$ (with the indexing convention $j+2n=j$) have the same
length. Then the quadrilateral $Q_j,1\le j\le n,$ built on any
pair of `opposite' sides of $P$ is inscribed into a circle, say
$C_j$. Let $\bo\Om$ be the closed curve formed by the appropriate
arcs in $C_j,1\le j\le n$. The domain $\Om=\Om(P)$ is obtained by
replacing each side of $P$ by an arc of the corresponding circle
$C_j$. Then $\Om$ is a floating domain with
$\de=1/2$~~\cite{ZaKo_french}. Moreover, assumption i) may be
replaced by the weaker assumption i') that in each of $n$ pairs of
the `opposite' sides of $P$ both sides have the same. The authors
state in~~\cite{ZaKo_french} that the class of polygons satisfying
assumptions i'), ii) is quite large. A simple example is the
$2n$-gon obtained from the regular $n$-gon by adding as vertices
the midpoints of its sides. Figure~~\ref{ZaKo_fig}, taken from the
russian original of~~\cite{ZaKo_french}, shows the floating domain
obtained by this construction from the equilateral
triangle.\footnote{The above construction of floating domains from
polygons is also contained in~~\cite{Ta06}, where it is ascribed
to N. Petrunin.}

\medskip

As for Question~~\ref{archimed_que} when $\de \ne 1/2$, there is a
controversy or, at least, confusion. According to
Tabachnikov~~\cite{Ta06}, E. Salkowski~~\cite{Sa34} claimed that
for $\de=\frac{m}{n}\ne\frac12$ the only convex floating domains
are the discs, but the proof is deficient.\footnote{The
work~~\cite{Sa34} is not available to the author.} According to
Auerbach~~\cite{Au38}, the ``last theorem'' in~~\cite{Sa34} states
that the only convex floating domains with
$\de=\frac{1}{n}\ne\frac12$ are the discs, but that he (Auerbach)
``unfortunately could not understand the proof''.  Furthermore,
Auerbach mentions that Salkowski separately proves the special
cases $\de=1/3,1/4$ of the ``last theorem'', and makes no comment
about the reliability of those proofs. On the other hand, Theorems
5 and 6 in~~\cite{Ta06} concern the cases $\de=1/3$ and $1/4$,
respectively, of Salkowski's ``last theorem''. The
work~~\cite{Ru_french} also contains a proof of the claim that
`any convex floating domain is a disc' in the case $\de=1/4$.
However, since~~\cite{Ru_french} proves this claim immediately
after ``proving'' the same for $\de=1/2$, the extremely sketchy
proof is likely to be deficient. There is also a considerable
literature on non-convex floating domains with $\de\ne1/2$
\cite{BMO04,Va09,We03}, but we will not review it here. The
interested reader may also consult the web site
$<$http://www.tphys.uni-heidelberg.de/wegner/Fl2mvs/Movies.html$>$
maintained by F. Wegner.

\medskip

\noindent{\bf Concluding remarks.} It seems appropriate to add a
few remarks of personal and social-historical character. I
acknowledge the anonymous referee of~~\cite{Gut12} who attracted
my attention to the literature on archimedean floating and
mentioned that equation~~\eqref{trig_eq} comes up in the
archimedean floating as well. See equations (14,18)
in~~\cite{Ta06}. This suggests the possibility of a nontrivial
relationship between the Finn-Young floating and the archimedean
floating. While visiting the Weizmann Institute of Science in
December 2010, I took part in the miniconference dedicated to the
$90$th birthday of professor Victor Abramovich Zalgaller. This
eventually led me to the russian originals of~~\cite{Ru_french}
and~~\cite{ZaKo_french}. The french
translations~~\cite{Ru_french,ZaKo_french} omit most of the
equations, as well as the figures. The transliteration {\em
Zalgaller} is more faithful to the original than ``Salgaller'',
used in ~~\cite{ZaKo_french}. In the russian alphabet the letter
corresponding to ``Z'' precedes ``K'', thus the coauthors in
~~\cite{ZaKo_french} are listed in the alphabetical order.

From Zalgaller's email communications I learned that the
archimedean floating problem was popular among older mathematics
students of Leningrad University while he and his friend Piotr
Kostelianets were sophomores. The results of the two sophomore
friends, as well as those of A.N. Ruban\footnote{He was a fifth
year student.} were simultaneously submitted for publication in
the Proceedings (Doklady) of the Soviet Academy of Sciences by
A.A. Markov, the head of the Geometry Chair at the
University.\footnote{Andrey Andreyevich Markov Jr. (1903 - 1979),
the son of A. A. Markov of the ``Markov chains'' fame, and an
outstanding mathematician as well.}

Shortly after, in June 1941, the war with the Nazi Germany began.
All three authors of~~\cite{Ru_french,ZaKo_french} went to the
front.\footnote{Zalgaller and Kostelianets volunteered in July
1941.} Kostelianets did not come back from the war. Although Ruban
survived, he became an invalid, no longer able to do mathematics.
Zalgaller was the luckiest of the three. Although severely
wounded, he was able to recover; he continued to fight in the
front lines until the Victory Day. Zalgaller resumed his
mathematical career in Leningrad, and became a distinguished
geometer. He co/authored many significant publications, including
several well known books. Zalgaller emigrated to Israel in
1999.\footnote{See
$<$http://en.wikipedia.org/wiki/VictorZalgaller$>$ and the
book~~\cite{Zal} for more information about Zalgaller.}

\medskip

Herman Auerbach (1901-1942) was a Polish geometer who lived in
Lw\'ow, which was a distinguished mathematical center during the
20 some years between the two World Wars. He was an active
participant in the mathematical activities at the Scottish Caf\'e
activities~~\cite{Scot}. It is plausible that Auerbach heard
Question~~\ref{archimed_que} personally from Ulam. In the years
1980-1984 I was friends with Marc Kac who studied in Lw\'ow in the
nineteen thirties. I much regret that I have never talked with
Marc about the `golden years of Scottish Caf\'e'. Alas, now it is
one of famous ``missed opportunities''~~\cite{missed}. In
September 1939 the Second World War broke out; Poland got invaded
by the Nazi Germany from the West and by the Soviet Union from the
East. Auerbach remained in Lw\'ow which became a part of the
Soviet Ukraine. He even published in the leading polish
mathematical journal, where the papers were now required to have
ukrainian summaries~~\cite{Au40}. In 1941 german troops captured
Lw\'ow. Auerbach perished in the hands of Gestapo in 1942. In 1992
the Polish Mathematical Society published a tribute to Herman
Auerbach and his work~~\cite{Au92i,Au92ii,Au92iii}. Auerbach wrote
the manuscript of~~\cite{Au92iii} in captivity, shortly before
execution.

%
%

%
%

%
%

\medskip

\noindent{\bf Acknowledgements.} I thank professor Zalgaller of
Rehovot, Israel, who sent me a scan of his personal copy of the
original of \cite{ZaKo_french} and patiently answered my
questions. The historical and personal information pertaining to
~~\cite{ZaKo_french} and~~\cite{Ru_french} came from these
communications, as well as from~~\cite{Al91} and~~\cite{Zal}. I
thank professor Alexey Borisov of Izhevsk, Russia, who provided me
with the originals of \cite{Ru_french,ZaKo_french}. I thank
professor Robert Finn of Palo Alto, California, who pointed out to
me a connection between the results of~~\cite{Gut93} and the
capillary floating and encouraged me to publish these results. The
work was partially supported by MNiSzW grant N N201 384834.

\end{document}

%% file: finn_young_float.pstex_t
\begin{picture}(0,0)%
\includegraphics{finn_young_float.pstex}%
\end{picture}%
\setlength{\unitlength}{2171sp}%
\begingroup\makeatletter\ifx\SetFigFont\undefined%
\gdef\SetFigFont#1#2#3#4#5{%
  \reset@font\fontsize{#1}{#2pt}%
  \fontfamily{#3}\fontseries{#4}\fontshape{#5}%
  \selectfont}%
\fi\endgroup%
\begin{picture}(12549,9154)(739,-9148)
\put(3001,-3211){\makebox(0,0)[lb]{\smash{{\SetFigFont{11}{13.2}{\familydefault}{\mddefault}{\updefault}{\color[rgb]{0,0,0}$\pi-\ga$}%
}}}}
\put(9376,-3211){\makebox(0,0)[lb]{\smash{{\SetFigFont{11}{13.2}{\familydefault}{\mddefault}{\updefault}{\color[rgb]{0,0,0}$\pi-\ga$}%
}}}}
\put(8626,-4936){\makebox(0,0)[lb]{\smash{{\SetFigFont{11}{13.2}{\familydefault}{\mddefault}{\updefault}{\color[rgb]{0,0,0}$\al$}%
}}}}
\end{picture}%

%% file: arch_float.pstex_t
\begin{picture}(0,0)%
\includegraphics{arch_float.pstex}%
\end{picture}%
\setlength{\unitlength}{2526sp}%
\begingroup\makeatletter\ifx\SetFigFont\undefined%
\gdef\SetFigFont#1#2#3#4#5{%
  \reset@font\fontsize{#1}{#2pt}%
  \fontfamily{#3}\fontseries{#4}\fontshape{#5}%
  \selectfont}%
\fi\endgroup%
\begin{picture}(7103,4912)(2213,-5525)
\put(3901,-5461){\makebox(0,0)[lb]{\smash{{\SetFigFont{10}{12.0}{\familydefault}{\mddefault}{\updefault}{\color[rgb]{0,0,0}$P(s)$}%
}}}}
\put(6376,-3886){\makebox(0,0)[lb]{\smash{{\SetFigFont{10}{12.0}{\familydefault}{\mddefault}{\updefault}{\color[rgb]{0,0,0}$C(s)$}%
}}}}
\put(6676,-4636){\makebox(0,0)[lb]{\smash{{\SetFigFont{10}{12.0}{\familydefault}{\mddefault}{\updefault}{\color[rgb]{0,0,0}$\Om(s)$}%
}}}}
\put(8326,-4936){\makebox(0,0)[lb]{\smash{{\SetFigFont{10}{12.0}{\familydefault}{\mddefault}{\updefault}{\color[rgb]{0,0,0}$A(s)$}%
}}}}
\put(4801,-5236){\makebox(0,0)[lb]{\smash{{\SetFigFont{8}{9.6}{\familydefault}{\mddefault}{\updefault}{\color[rgb]{0,0,0}$\het(s)$}%
}}}}
\put(2926,-3061){\makebox(0,0)[lb]{\smash{{\SetFigFont{12}{14.4}{\familydefault}{\mddefault}{\updefault}{\color[rgb]{0,0,0}$\Om$}%
}}}}
\put(5326,-736){\makebox(0,0)[lb]{\smash{{\SetFigFont{12}{14.4}{\familydefault}{\mddefault}{\updefault}{\color[rgb]{0,0,0}$\bo\Om$}%
}}}}
\put(9301,-2761){\makebox(0,0)[lb]{\smash{{\SetFigFont{10}{12.0}{\familydefault}{\mddefault}{\updefault}{\color[rgb]{0,0,0}$P(s+\de)$}%
}}}}
\put(8476,-3361){\makebox(0,0)[lb]{\smash{{\SetFigFont{8}{9.6}{\familydefault}{\mddefault}{\updefault}{\color[rgb]{0,0,0}$\het(s)$}%
}}}}
\end{picture}%

%% file: za_ko.pstex_t
\begin{picture}(0,0)%
\includegraphics{za_ko.pstex}%
\end{picture}%
\setlength{\unitlength}{2526sp}%
\begingroup\makeatletter\ifx\SetFigFont\undefined%
\gdef\SetFigFont#1#2#3#4#5{%
  \reset@font\fontsize{#1}{#2pt}%
  \fontfamily{#3}\fontseries{#4}\fontshape{#5}%
  \selectfont}%
\fi\endgroup%
\begin{picture}(7230,6395)(2386,-8291)
\put(9601,-2386){\makebox(0,0)[lb]{\smash{{\SetFigFont{10}{12.0}{\familydefault}{\mddefault}{\updefault}{\color[rgb]{0,0,0}$A$}%
}}}}
\put(2401,-2461){\makebox(0,0)[lb]{\smash{{\SetFigFont{10}{12.0}{\familydefault}{\mddefault}{\updefault}{\color[rgb]{0,0,0}$B$}%
}}}}
\put(5926,-8236){\makebox(0,0)[lb]{\smash{{\SetFigFont{10}{12.0}{\familydefault}{\mddefault}{\updefault}{\color[rgb]{0,0,0}$C$}%
}}}}
\put(5926,-2161){\makebox(0,0)[lb]{\smash{{\SetFigFont{10}{12.0}{\familydefault}{\mddefault}{\updefault}{\color[rgb]{0,0,0}$C'$}%
}}}}
\put(4126,-5461){\makebox(0,0)[lb]{\smash{{\SetFigFont{10}{12.0}{\familydefault}{\mddefault}{\updefault}{\color[rgb]{0,0,0}$A'$}%
}}}}
\put(7876,-5461){\makebox(0,0)[lb]{\smash{{\SetFigFont{10}{12.0}{\familydefault}{\mddefault}{\updefault}{\color[rgb]{0,0,0}$B'$}%
}}}}
\end{picture}%